\documentclass[12pt,titlepage]{article}
\usepackage{amssymb}
\usepackage{graphicx}
\usepackage{caption2}

\newcommand{\be}{\begin{equation}}
\newcommand{\ee}{\end{equation}}
\newcommand{\bea}{\begin{eqnarray}}
\newcommand{\eea}{\end{eqnarray}}
\newcommand{\bef}{\begin{figure}}
\newcommand{\ef}{\end{figure}}
\newcommand{\bt}{\begin{tabular}}
\newcommand{\et}{\end{tabular}}
\newcommand{\bno}{\begin{enumerate}}
\newcommand{\eno}{\end{enumerate}}

\def\3{\ss}

\catcode`\"=\active
\def"{\accent'177}

\begin{document}
\markright{\hspace{2in} G. Shabbir, H. Khan and M. A. Sadiq}
\thispagestyle{plain}
\begin{center}
\bf \Large{A note on Exact solution of SIR and SIS epidemic
models}
\end{center}

\begin{center}
{ G. Shabbir$^1$\footnote{ E-mail address: shabbir@giki.edu.pk },
H. Khan$^1$ \footnote{The corresponding author, E-mail address:
{hina.maths@gmail.com}}  and M. A. Sadiq $^2$ \footnote{ E-mail
address: {msa502@york.ac.uk}}
  \hspace{.2in}\/}

\small $^1$ Faculty of Engineering Sciences,\\
GIK Institute of Engineering Sciences and Technology,\\
 Topi, Swabi, Khyber Pakhtunkhwa Pakistan.\\
\small $^2$ Department of Mathematics, University of York, UK.
\end{center}

{\bf Abstract}

{\em In this article we have successfully obtained an exact
solution of a particular case of SIR and SIS epidemic models given
by Kermack and Mckendrick \cite{Kermack-1927} for constant
population, which are described by coupled nonlinear differential
equations. Our result has no limiting conditions for any parameter
involved in the given models. In epidemiology many researchers
believe that it is very hard to get an exact solution for such
models. We hope this solution will be an opening window and good
addition in the area of epidemiology.
}\\


\section{Introduction}
The mathematical modelling of spread of diseases is an area of
biology named as epidemiology. Physical real life problems arising
in epidemiology may be described, in a first formulation, using
differential equations. As far as our knowledge is concerned the
first mathematical model of epidemiology was formulated and solved
by Daniel Bernoulli in 1760. Since the time of Kermack and
McKendrick \cite{Kermack-1927}, the study of mathematical
epidemiology has grown rapidly, with a large variety of models
having been formulated and applied to infectious diseases
\cite{Diekmann-2007}-\cite{West-1997}.\\

 The dilemma with models in epidemiology is,
sometimes even the simplest mathematical model of natural
phenomena with sets of first-order ordinary differential equations
are non-integrable. Nucci and Leach \cite{nucci-2004} have
discussed these features and obtained an integrable SIS model by
applying Lie analysis. According to Nucci and Leach
\cite{nucci-2004} {\em It would seem that a fatal disease which
this models is also not good for mathematics!}. Khan et.al.,
\cite{khan-2009} has obtained the convergent series solutions for
SIR and SIS models using homotopy analysis method. Considerable
attention has been given to these models by several authors using
stability, bifurcation theory, Lyapunov method,
Poincar$\acute{e}$-Bendixson type theorems, index and topological
concepts see (\cite{Diekmann-2007}-\cite{Murray1992}).

To understand what SIS and SIR stands for let us consider some
constant population, let divide them into three components: the
susceptible component, S, who can catch the disease; the infective
component, I, who are infected and can transmit the disease to
susceptible, and the removed component, R, who either had the
disease and recovered, died or have developed immunity, or have
been removed from contact with other components. SIR model
contains all three classes, while, in SIS model the infectives
return to the susceptible class on recovery because the disease
confers no immunity against reinfection. Such models are more
effective for disease caused by bacteria or
helminth agents and also for most sexually transmitted disease.\\

In fact closed-form solutions to mathematical models play a vital
role in the proper understanding of qualitative features of
natural science. In epidemiology mathematical models describes
some effects of disease in some population so one need to be very
careful while analyzing such models. Keeping this view of point,
here in this article our main goal is to obtain pure analytic
exact solution to the SIS and SIR epidemic model which is a
particular case of the model given by Kermack and McKendrick
\cite{Kermack-1927}.

In second section we present the solution for SIS model while in
third section we present solution for SIR model.

\section{Mathematical Analysis for SIS model}
 The SIS model by Kermack and Mckendrick \cite{Kermack-1927} is given as follows
\begin{equation}
s'(t)=-rsi+{\alpha}i,\label{ffirst:m}
\end{equation}
\begin{equation}
i'(t)=rsi-{\alpha}i,\label{22nd:m}
\end{equation}
subject to Initial conditions,
\begin{equation}
i(0)=I_0,\;\;\;s(0)=S_0,\label{bBC1:m}
\end{equation}
where $r>0$, $I_0>0$ and $S_0>0$. Here $r$ is the infectivity
coefficient of the typical Lotka-Volterra interaction term and
$\alpha$ the recovery coefficient. From equation (\ref{ffirst:m})
and (\ref{22nd:m}) one can easily observe that $s(t)+i(t)=k$,
where $k$ is total population. This model is different from SIR
model as in this model the recovered members return to the class S
at a rate of $i\alpha$
instead of moving to class R.\\
Now using $s(t)=k-i(t)$ in equation (\ref{22nd:m}),

\begin{eqnarray}
{i'}=r\;(k-i)i-\alpha i.
 \label{ddef:N}
\end{eqnarray}
Integrating the above equation (\ref{ddef:N}) with substitution
$y=i^{-1}$, we get

\begin{eqnarray}
 ye^{(rk-\alpha)t}=\int{re^{(rk-\alpha)t}dt}+C,\label{ddef:N2}
\end{eqnarray}

where $C$ is integrating constant. After simplifying equation.
(\ref{ddef:N2}) and letting $\beta =rk-\alpha$, we get,
\begin{eqnarray}
y=\frac{r+\beta C e^{-\beta t}}{\beta}. \label{def:nN23}
\end{eqnarray}
Using back substitution $y=i^{-1}$,
\begin{eqnarray}
i(t)=\frac{\beta}{r+\beta C e^{-\beta t}}. \label{def:nN233}
\end{eqnarray}
After applying initial condition $i(0)=i_0$ on equation
(\ref{def:nN233}) we get,
\begin{eqnarray}
C=\frac{\beta-i_0r}{\beta i_0}\nonumber. \label{def:nN234}
\end{eqnarray}
The exact solution for $s(t)$ obtained  by using the fact
$s(t)+i(t)=k$ is,
\begin{eqnarray}
s(t)=k-\frac{\beta}{r+\beta (\frac{\beta-i_0r}{\beta i_0})
e^{-\beta t}}\nonumber. \label{def:nN23r3}
\end{eqnarray}
\section{Mathematical Analysis for SIR model}
 The SIR model by Kermack and Mckendrick \cite{Kermack-1927} is given as follows
\begin{equation}
s'(t)=-\beta s(t)i(t)-\mu s(t)+{\mu},\label{first:m}
\end{equation}
\begin{equation}
i'(t)=\beta s(t)i(t)-{\mu}i(t),\label{2nd:m}
\end{equation}
subject to Initial conditions,
\begin{equation}
i(0)=i_0,\;\;\;s(0)=s_0,\label{BC1:m}
\end{equation}
$(')$ denotes the derivative with respect to time, where
$\beta>0$, $i_0>0$ and $s_0>0$. Here $\beta$ is the infectivity
coefficient of the typical Lotka-Volterra interaction term and
$\mu$ the recovery coefficient.  Now adding equation
(\ref{first:m}) and equation (\ref{2nd:m}),

\begin{eqnarray}
{(s+i)'}=-\mu\;(s+i)+\mu.
 \label{def:N}
\end{eqnarray}
Integrating the above equation \ref{def:N}, we get
\begin{eqnarray}
s(t)=1+Ce^{(-\mu)t}-i(t),\label{def:N2}
\end{eqnarray}

where $C$ is an integrating constant. Now using equation
(\ref{def:N2}) in equation (\ref{2nd:m}), we get,
\begin{eqnarray}
i'(t)=(\beta -\mu+\beta C e^{-\mu t})i(t)-\beta i^2(t).
\label{def:N23}
\end{eqnarray}
Integrating the above equation with substitution $z=i^{-1}$, we
get
\begin{eqnarray}
 ze^{(\beta-\mu)t-(\frac{\beta C}{\mu})e^{-\mu t}}=
 \int{\beta e^{(\beta-\mu)t-(\frac{\beta C}{\mu})e^{-\mu t}}dt}+D.\label{def:N72}
\end{eqnarray}
Where $D$ is an integrating constant. Now observing that $e^{- \mu
t}$ approaches to zero as $t\rightarrow\infty$, here we take
series expansion for $e^{- \mu t}$ and neglecting square and
higher terms we get $e^{- \mu t}=1-\mu t$. Substituting back in
equation (\ref{def:N72}) one gets,
\begin{eqnarray}
 ze^{(\beta-\mu)t-(\frac{\beta C}{\mu})(1-\mu t)}=
 \int{\beta e^{(\beta-\mu)t-(\frac{\beta C}{\mu})(1-\mu t)}dt}+D,\label{def:N732}
\end{eqnarray}
\begin{eqnarray}
 z=\frac{\beta}{\beta-\mu+\beta C}+D e^{-(\beta-\mu+\beta C)t}e^{\frac{\beta C}{\mu}}
 ,\label{def:NN732}
\end{eqnarray}
Substituting $z=i^{-1}$ and $\beta-\mu +\beta C= \lambda$, in
equation (\ref{def:NN732}) we get
\begin{eqnarray}
i(t)=\frac{\lambda}{\beta +\lambda D e^{-\lambda t}e^{\frac{\beta
C}{\mu}}}. \label{def:N233}
\end{eqnarray}
Applying initial condition $i(0)=i_0$ and $s(0)=s_0$ in equation
(\ref{def:N2}) and (\ref{def:N233}) we get,
\begin{eqnarray}
D=\frac{\lambda-i_0 \beta}{\lambda i_0e^{\frac{\beta
C}{\mu}}}\nonumber\\C=s_0+i_0-1. \label{def:N234}
\end{eqnarray}
The exact solution for $s(t)$ obtained  by using the value of
$i(t)$ in equation (\ref{def:N2}) which is,
\begin{eqnarray}
s(t)=1+(s_0+i_0-1 ){(1-\mu t)}-\frac{\lambda}{\beta +\lambda
(\frac{\lambda-i_0 \beta}{\lambda i_0e^{\frac{\beta
(s_0+i_0-1)}{\mu}}}) e^{-\lambda t+\frac{\beta
(s_0+i_0-1)}{\mu}}},\label{def:N122}
\end{eqnarray}
\begin{eqnarray}
i(t)=\frac{\lambda}{\beta +\lambda (\frac{\lambda-i_0
\beta}{\lambda i_0e^{\frac{\beta (s_0+i_0-1)}{\mu}}}) e^{-\lambda
t+\frac{\beta (s_0+i_0-1)}{\mu}}}. \label{def:N233}
\end{eqnarray}
\section{Conclusion}In this paper an exact analytic solution to a
particular case of SIR and SIS models with constant population is
obtained with the help of direct integration tools. As these
models have already been observed by different researchers with
different techniques but we know the importance of exact solution,
that exact solutions allow researchers to run experiments, by
introducing natural initial or boundary conditions, to observe the
behavior of functions involved in equations. This present work
will be of great support to the epidemiology to investigate the
hidden phenomenon and validity of these models which is already
being considered and used by different researchers. This kind of
approach will
attract other nonlinear models as well having been analyzed again for the possible existence of exact solution.  \\

{\bf Acknowledgement}\\

Many thanks to Prof. Dr. S. J. Liao for bringing the attentions towards the area of epidemiology.

\end{document}